\documentclass[12pt]{article}
\usepackage[T2A]{fontenc}
\usepackage[utf8]{inputenc}
\usepackage{amssymb,amsmath,amsfonts,amsthm,amscd,latexsym,verbatim,graphics,epsfig,indentfirst,color}
\usepackage{geometry}
\def\id{\mathrm{id}}
\def\rb{\textrm{rb}}

\def\rad{\mathrm{rad}}

\def\Span{\mathrm{Span}}

\def\Rad{\mathrm{Rad}}
\def\Aut{\mathrm{Aut}}

\def\Imm{\mathrm{Im}}

\begin{document}

\sloppy

\begin{center}
{\Large
Rota---Baxter operators of weight zero \\
on the matrix algebra of order three without unit in kernel}

Vsevolod Gubarev
\end{center}

\begin{abstract}
We describe all Rota---Baxter operators~$R$ of weight zero
on the matrix algebra~$M_3(F)$ over a quadratically closed field~$F$ of characteristic not 2 or 3 such that $R(1)\neq0$.
Thus, we get a partial classification of solutions to the associative Yang---Baxter equation on~$M_3(F)$. 
For the solution, the computer algebra system \texttt{Singular} was involved.

\medskip


\end{abstract}

\section{Introduction}

Given an algebra $A$ and a scalar $\lambda\in F$, where $F$ is a~ground field,
a~linear operator $R\colon A\rightarrow A$ is called a Rota---Baxter operator
(RB-operator, for short) on $A$ of weight~$\lambda$ if the following identity
\begin{equation}\label{RB}
R(x)R(y) = R( R(x)y + xR(y) + \lambda xy )
\end{equation}
holds for all $x,y\in A$. Then the algebra $A$ is called a~Rota---Baxter algebra (RB-algebra).

Glen Baxter in 1960~\cite{Baxter} introduced the notion of a Rota---Baxter operator
as a~formal generalization of the integration by parts formula.
Further, J.B.~Miller~\cite{Miller66}, G.-C.~Rota~\cite{Rota}, P.~Cartier~\cite{Cartier} and others studied such operators on commutative algebras.

In the 1980s, the deep connection between solutions of the classical Yang---Baxter equation from mathematical physics and RB-operators of weight zero on a semisimple finite-dimensional Lie algebra was found by A.A. Belavin and V.G. Drinfel'd~\cite{BelaDrin82} and M.A.~Semenov-Tyan-Shanskii~\cite{Semenov83}.
Actually, they rediscovered Rota---Baxter operators via the classical Yang---Baxter equation.
The same correspondence applied for solutions of the modified Yang---Baxter equation and Rota---Baxter operators of nonzero weight on a simple finite-dimensional Lie algebra also works~\cite{Reshetikhin,Gonch}.

From the 2000s, the active study of Rota---Baxter operators on associative algebras has  begun; see the monograph~\cite{GuoMonograph}.
The main example of the study is the matrix algebra.
For the algebra $M_n(\mathbb{C})$, solutions of the associative Yang---Baxter equation~\cite{Aguiar2000b,Zhelyabin} are in one-to-one correspondence with RB-operators of weight zero~\cite{Unital}.
Also, skew-symmetric Rota---Baxter operators of weight zero on~$M_n(\mathbb{C})$
are in one-to-one correspondence with finite-dimensional double Lie algebras~\cite{DoubleLie}
arisen in noncommutative geometry.

Thus, the problem of the classification of Rota---Baxter operators on the matrix algebra is a very natural one.
They were classified on $M_2(\mathbb{C})$ \cite{BGP,Mat2}.
In 2013, V.V.~Sokolov described all skew-symmetric RB-operators of nonzero weight on $M_3(\mathbb{C})$~\cite{Sokolov}; up to conjugation with automorphisms and transpose, he got 8 series.
In the series of the three works~\cite{GonGub,Gub2021,Gub2024},
the complete classification of Rota---Baxter operators of nonzero weight~$\lambda$ on $M_3(F)$ was obtained.
The current work is devoted to the description of RB-operators~$R$ of weight~0 on $M_3(F)$ such that $R(1) \neq 0$.
Actually, we get 24 series, some of them involve parameters.

In~\cite{Unital}, some general properties of RB-operators were stated. In particular, for every RB-operator $R$ of weight zero on $M_n(F)$, we have $R^{2n-1} = 0$ and $R(1)$ is nilpotent.
In~\cite{Spectrum}, this fact was generalized to any weight as follows.
Given a~Rota---Baxter operator $R$ of weight~$\lambda$ on $M_n(F)$, we have $R^k(R+\lambda\id)^l = 0$
for some suitable $k,l$. The property concerned the minimal polynomial of $R(1)$, where $R$ is an RB-operator of nonzero weight~$\lambda$ on~$M_n(F)$, was applied in~\cite{GonGub} to describe all RB-operators of nonzero weight~$\lambda$ on~$M_n(F)$, which are not projections on a subalgebra parallel to another one.
The remaining case of projective RB-operators of nonzero weight on~$M_n(F)$ was considered in~\cite{Gub2021,Gub2024}.

Returning to the case of weight zero and $n = 3$, the main tool for us is that $R(1)$ is nilpotent.
Up to conjugation with an automorphism of~$M_3(F)$, we may assume that 
\begin{itemize}
 \item $R(1) = 0$, 
 \item $R(1) = e_{12}$, 
 \item $R(1) = e_{12} + e_{13}$,
\end{itemize}
Here $e_{ij}$ denotes the matrix unity. 

In the current work, we study the last two cases. We put $R$ to be a linear operator with unknown coefficients. After some initial steps, we study the system of quadratic equations on these coefficients.
For this, we actively use the computer algebra system \texttt{Singular}~\cite{DGPS}. 

The work on the first case ($R(1) = 0$) is in progress.

Let us give a brief outline of the work.
In~\S2, we provide the required preliminaries.
In~\S3, we study the case $R(1) = e_{12} + e_{13}$ and show that up to conjugation with automorphisms $M_2(F)$, transpose,
and multiplication on a nonzero scalar, there is the only one appropriate Rota---Baxter operator~(Q1). 
For this, we apply the result from the work~\cite{Upper-triangular}, where Rota---Baxter operators of weight~0 on the subalgebra $U_3(F)$ of upper-triangular matrices were described.

In~\S4, we deal with the case $R(1) = e_{12}$.
Firstly, we show that we may assume that $e_{12}, e_{13} \in \ker R$, 
$R(e_{32}),R(e_{11}),R(e_{22})\in \Span\{e_{12},e_{13}\}$, $\Imm R\subset \Span\{e_{11},e_{1i},e_{2i},e_{3i}\mid i=2,3\}$.
Secondly, we consider different cases depending on the values of the coefficients of $\Pr(R(e_{32}))|_{e_{12}}$, $\Pr(R(e_{32}))|_{e_{13}}$ and so on.
Here we get 23 operators~(Q2)--(Q24), some of them involve two-valued parameters and at most one parameter taking values from~$F$.
Note that we intensively use conjugation with automorphisms of~$M_3(F)$ preserving~$e_{12}$.
For this, we completely describe such (anti)automorphisms of~$M_3(F)$. 
We collect all obtained Rota---Baxter operators~(Q1)--(Q24) in Theorem~3.
Finally, we show that all RB-operators from Sokolov's classification of skew-symmetric RB-operators on~$M_3(F)$~\cite{Sokolov} such that $R(1)\neq0$ (there are exactly seven of these out of eight) are among the operators from Theorem~3.

We assume that the ground field~$F$ is of characteristic not 2 or 3 and is quadratically closed.

\section{Preliminaries}

We write down some basic properties of Rota---Baxter operators.

{\bf Proposition 1}~\cite{GuoMonograph}.
Let $A$ be an algebra, and let $R$ be an RB-operator of weight~0 on~$A$.
Then the operator $\alpha R$ for any $\alpha\in F$ is again an RB-operator of weight~0 on~$A$.

{\bf Proposition 2}~\cite{GuoMonograph}.
Let $A$ be an algebra, and let $R$ be an RB-operator of weight~0 on~$A$.
Then $\Imm(R)$ is a subalgebra of~$A$ and $\ker R$ is an $\Imm(R)$-bimodule. 

{\bf Proposition 3}~\cite{BGP}.
Given an algebra $A$, an RB-operator $R$ of weight $\lambda$ on $A$,
and $\psi\in\Aut(A)$, the operator $R^{(\psi)} = \psi^{-1}R\psi$
is again an RB-operator of weight $\lambda$ on $A$.

The same result is true when $\psi$ is an antiautomorphism of $A$, which means an automorphism of the vector space~$(A,+)$
such that $\psi(xy) = \psi(y)\psi(x)$ for all $x,y\in A$.
Transpose on the matrix algebra gives an example of an antiautomorphism of $M_n(F)$.

{\bf Proposition 4}~\cite{BGP,GuoMonograph}.
Let $A$ be a unital algebra, and let $R$ be an RB-operator of weight~0 on~$A$.

(a) If $R(1) = 0$, then $\Imm R\subset \ker R$.

(b) We have $(R(1))^n = n! R^n(1)$, $n\in\mathbb{N}$.

Now, we put some essential properties of RB-operators.

{\bf Lemma 1}~\cite{Unital}.
Let $A$ be a unital associative finite-dimensional algebra and $R$ be an RB-operator
on $A$ of weight zero. Then $R(1)$ is nilpotent.

{\bf Theorem 1}~\cite{Unital,Spectrum,Miller66}.
Let $A$ be a unital associative finite-dimensional algebra.
Then there exists $N$ such that $R^N = 0$
for every RB-operator~$R$ of weight~0 on $A$.

Let $A$ be an algebra.
Due to~\cite{Unital}, a Rota---Baxter index $\rb(A)$ of $A$ is defined as follows,
$$
\rb(A) = \min\{n\in\mathbb{N}\mid R^n = 0\mbox{ for an RB-operator }R \mbox{ of weight zero on }A\}.
$$
If such a number is undefined, put $\rb(A) = \infty$.

By Theorem~1, $\rb(M_n(F))<\infty$.
In~\cite{Unital}, it was shown that $\rb(M_n(F)) = 2n + 1$.
The next example provides the lower bound.

{\bf Example 1}~\cite{Unital}.
A linear operator $R$ defined on $M_n(F)$ as follows
\begin{equation}\label{MaxRbMat}
R(e_{ij}) = \begin{cases}
e_{i,j+1}+e_{i+1,j+2}+\ldots+e_{n-j+i-1,n}, & i\leq j\leq n-1, \\
-(e_{i-1,j}+e_{i-2,j-1}+\ldots+e_{i-j+1,1}), & i>j, \\
0, & j = n, \end{cases}
\end{equation}
is an RB-operator of weight zero.
Due to the definition, we have $R^{2n-1} = 0$ and $R^{2n-2}\neq0$.

Below, we write down the classification of Rota---Baxter operators of weight~0 on $M_2(F)$.

{\bf Theorem 2}~\cite{Aguiar2000b,BGP,Mat2,Panasenko}.
All nonzero RB-operators of weight zero on $M_2(F)$ over a~field~$F$ 
of characteristic not two up to conjugation with automorphisms $M_2(F)$, transpose
and multiplication on a nonzero scalar are the following:

(L1) $R(e_{21}) = e_{12}$, $R(e_{11}) = R(e_{12}) = R(e_{22}) = 0$;

(L2) $R(e_{21}) = e_{11}$, $R(e_{11}) = R(e_{12}) = R(e_{22}) = 0$;

(L3) $R(e_{21}) = e_{11}$, $R(e_{22}) = e_{12}$, $R(e_{11}) = R(e_{12}) = 0$;

(L4) $R(e_{21}) =- e_{11}$, $R(e_{11}) = e_{12}$, $R(e_{12}) = R(e_{22}) = 0$.

\noindent Moreover, these operators lie in different orbits
of the set of RB-operators of weight~0 on $M_2(F)$ under
conjugation with $\Aut(M_2(F))$, transpose or multiplication on a nonzero scalar.

Below, we will apply the following automorphism $\psi_{r,s}$ of $M_3(F)$ for $r,s\in F\setminus \{ 0 \}$:
\begin{equation} \label{psi}
\begin{gathered}
\psi_r(e_{ii}) = e_{ii}, \ i=1,2,3,\quad
\psi_r(e_{12}) = re_{12}, \quad
\psi_r(e_{23}) = se_{23}, \\
\psi_r(e_{13}) = rs e_{13}, \quad
\psi_r(e_{21}) = (1/r)e_{21}, \quad
\psi_r(e_{32}) = (1/s)e_{32}, \quad
\psi_r(e_{31}) = (1/rs)e_{31}.
\end{gathered}
\end{equation}

\section{RB-operators on $M_3(F)$ of maximal rb-index}

Consider the case
\begin{equation}\label{UnitImage}
R(1) = e_{12} + e_{23}.
\end{equation}

By~\cite[Lemma\,5.23]{Unital}, the subalgebra of upper-triangular matrices $U_3(F)$ is $R$-inva- \linebreak
riant.
By~\cite{Upper-triangular}, $R$ is conjugate with the help of $\Aut(U_3(F))$ to the following operator:
\begin{equation} \label{R1=e12+e23-FirstForm}
\begin{gathered}
R(e_{13}) = R(e_{23}) = 0,\quad
R(e_{12}) = \frac{1}{2}e_{13},\\
R(e_{11}) = \begin{pmatrix}
0 & 1 & r_{13} \\
0 & 0 & 1/2 \\
0 & 0 & 0 \\
\end{pmatrix},\
R(e_{22}) = \begin{pmatrix}
0 & 0 & p_{13} \\
0 & 0 & 1/2 \\
0 & 0 & 0 \\
\end{pmatrix},\
R(e_{33}) = \begin{pmatrix}
0 & 0 & -r_{13}-p_{13} \\
0 & 0 & 0 \\
0 & 0 & 0 \\
\end{pmatrix}.
\end{gathered}
\end{equation}
Since all automorphisms of $U_3(F)$ are inner, we may apply this result.

Applying RB-identity, we have
\begin{gather}
R(x)(e_{12}+e_{23})
= R(x)R(1) = R^2(x) + R(x(e_{12}+e_{23})), \label{RxR1} \\
(e_{12}+e_{23})R(x)
= R(1)R(x) = R^2(x) + R((e_{12}+e_{23})x). \label{R1Rx}
\end{gather}
Denote $A = (a_{ij}) = R(x)$, then
\begin{equation}\label{R^2x}
R^2(x) = \begin{pmatrix}
0 & a_{11} & a_{12} \\
0 & a_{21} & a_{22} \\
0 & a_{31} & a_{32} \\
\end{pmatrix} - R(x(e_{12}+e_{23}))
= \begin{pmatrix}
a_{21} & a_{22} & a_{23} \\
a_{31} & a_{32} & a_{33} \\
0 & 0 & 0 \\
\end{pmatrix} - R((e_{12}+e_{23})x).
\end{equation}

Put $x = e_{21}$ in \eqref{R^2x}: $R(e_{21})R(1) - R(e_{22}) = R(1)R(e_{21}) - R(e_{11})$, so
one has for $A = R(e_{21})$ the following equalities:
$$
a_{21} = a_{31} = a_{32} = 0,\quad
a_{22} = a_{33},\quad
a_{22} = a_{11} + 1,\quad
a_{12} = a_{23} - r_{13} + p_{13}.
$$
Since $\Imm(R)$ contains only degenerate matrices, we have
$$
R(e_{21}) = \begin{pmatrix}
a_{22}-1 & a_{23}-r_{13}+p_{13} & a_{13} \\
0 & a_{22} & a_{23} \\
0 & 0 & a_{22} \\
\end{pmatrix},\quad a_{22}\in\{0,1\}.
$$

Further, consider
\begin{multline*}
R(e_{21})R(e_{21}) \\
 = \begin{pmatrix}
(a_{22}-1)^2 & (2a_{22}-1)(a_{23}-r_{13}+p_{13}) & (2a_{22}-1)a_{13}+a_{23}(a_{23}-r_{13}+p_{13}) \\
0 & a_{22}^2 & 2a_{22}a_{23} \\
0 & 0 & a_{22}^2 \\
\end{pmatrix};
\end{multline*}
\begin{multline*}
R(R(e_{21})e_{21} + e_{21}R(e_{21})) \\
 = R( (a_{23}-r_{13}+p_{13})e_{11} + a_{22}e_{21} + (a_{22}-1)e_{21} + (a_{23}-r_{13}+p_{13})e_{22} + a_{13}e_{23}) \\
 = (2a_{22}-1)R(e_{21}) + (a_{23}-r_{13}+p_{13})R(e_{11}+e_{22}).
\end{multline*}
Comparing the coefficients on matrix unities, we get
\begin{gather*}
2a_{22}a_{23} = (2a_{22}-1)a_{23} + (a_{23}-r_{13}+p_{13}), \\
(2a_{22}-1)(a_{23}-r_{13}+p_{13})
 = (2a_{22}-1)(a_{23}-r_{13}+p_{13}) + (a_{23}-r_{13}+p_{13}),
\end{gather*}
so $r_{13} = p_{13}$ and $a_{23} = 0$.

Put $x = e_{32}$ in \eqref{R^2x}: $R(e_{32})R(1) - R(e_{33}) = R(1)R(e_{32}) - R(e_{22})$, 
thus, for $B = R(e_{32})$ we have
$$
b_{21} = b_{31} = b_{32} = 0,\quad
b_{11} = b_{22},\quad
b_{33} = b_{22}+1/2,\quad
b_{12} = b_{23} - 3r_{13}.
$$
Since $\Imm(R)$ contains only degenerate matrices, we have
$$
R(e_{32})
 = \begin{pmatrix}
b_{11} & b_{23}-3r_{13} & b_{13} \\
0 & b_{11} & b_{23} \\
0 & 0 & b_{11}+1/2 \\
\end{pmatrix},\quad b_{11}\in\{0,-1/2\}.
$$

Let us check the RB-identity for $x = y = e_{32}$:
$$
R(e_{32})R(e_{32})
 = \begin{pmatrix}
b_{11}^2 & 2b_{11}(b_{23}-3r_{13}) & (2b_{11}+1/2)b_{13}+b_{23}(b_{23}-3r_{13}) \\
0 & b_{11}^2 & (2b_{11}+1/2)b_{23} \\
0 & 0 & (b_{11}+1/2)^2 \\
\end{pmatrix};
$$

\vspace{-0.3cm}

\begin{multline*}
R(R(e_{32})e_{32} + e_{32}R(e_{32})) 
 = R( b_{13}e_{12} + b_{23}e_{22} + (b_{11}+1/2)e_{32} + b_{11}e_{32} + b_{23}e_{33} ) \\
 = (2b_{11}+1/2)R(e_{32}) + \frac{b_{13}}{2}e_{13} + b_{23}R(e_{22}+e_{33}).
\end{multline*}
Comparing the coefficients on matrix unities, we get
\begin{gather*}
2b_{11}(b_{23}-3r_{13}) = (2b_{11}+1/2)(b_{23}-3r_{13}), \\
(2b_{11}+1/2)b_{23}
 = (2b_{11}+1/2)b_{23} + b_{23}/2, \\
(2b_{11}+1/2)b_{13}+b_{23}(b_{23}-3r_{13})
 = (2b_{11}+1/2)b_{13} + b_{13}/2 - r_{13} b_{23},
\end{gather*}
so $b_{23} = 3r_{13} = 0$ and $b_{13} = 0$.

Hence,
\begin{gather*}
R(e_{13}) = R(e_{23}) = R(e_{33}) = 0,\quad
R(e_{12}) = \frac{1}{2}e_{13},\quad
R(e_{11}) = e_{12} + \frac{1}{2}e_{23},\quad
R(e_{22}) = \frac{1}{2}e_{23},\\
R(e_{21})
 = \begin{pmatrix}
a_{22}-1 & 0 & a_{13} \\
0 & a_{22} & 0 \\
0 & 0 & a_{22} \\
\end{pmatrix},\quad
R(e_{32})
 = \begin{pmatrix}
b_{11} & 0 & 0 \\
0 & b_{11} & 0 \\
0 & 0 & b_{11}+1/2 \\
\end{pmatrix},
\end{gather*}
where $a_{22}\in\{0,1\}$ and $b_{11}\in\{0,-1/2\}$.

Finally, put $x = e_{31}$ in \eqref{R^2x}:
$R(e_{31})R(1) - R(e_{32}) = R(1)R(e_{31}) - R(e_{21})$,
so for $C = R(e_{31})$ we have the matrix equality
$$
\begin{pmatrix}
0 & c_{11} & c_{12} \\
0 & c_{21} & c_{22} \\
0 & c_{31} & c_{32} \\
\end{pmatrix}
 + \begin{pmatrix}
a_{22}-b_{11}-1 & 0 & a_{13} \\
0 & a_{22}-b_{11} & 0 \\
0 & 0 & a_{22}-b_{11}-1/2 \\
\end{pmatrix}
 = \begin{pmatrix}
c_{21} & c_{22} & c_{23} \\
c_{31} & c_{32} & c_{33} \\
0 & 0 & 0 \\
\end{pmatrix},
$$
which implies the following relations
\begin{gather*}
c_{31} = 0,\quad
c_{21} = a_{22}-b_{11}-1,\quad
c_{11} = c_{22} = c_{33},\\
c_{23} = c_{12} + a_{13},\quad
c_{32} = c_{21} + a_{22}-b_{11},\quad
c_{32} = b_{11}+1/2-a_{22}.
\end{gather*}
In particular,
$c_{21} = a_{22}-b_{11}-1
 = 2b_{11}-2a_{22} +1/2$, i.\,e.,
$b_{11}-a_{22}+\frac{1}{2} = 0$. Since $a_{22}\in\{0,1\}$ and $b_{11}\in\{0,-1/2\}$,
the only possibility is $b_{11}-a_{22}+\frac{1}{2} = 0$ and $a_{22} = 0$, $b_{11} = -1/2$.

Therefore,
$$
R(e_{31})
 = \begin{pmatrix}
c_{11} & c_{23}-a_{13} & c_{13} \\
-1/2 & c_{11} & c_{23} \\
0 & 0 & c_{11} \\
\end{pmatrix}.
$$

Write down
\begin{multline*}
R(e_{31})R(e_{31})
= \begin{pmatrix}
c_{11}^2-(c_{23}-a_{13})/2 & 2c_{11}(c_{23}-a_{13}) & 2c_{11}c_{13} + c_{23}(c_{23}-a_{13}) \\
-c_{11} & c_{11}^2 - (c_{23}-a_{23})/2 & 2c_{11}c_{23} - c_{13}/2 \\
0 & 0 & c_{11}^2 \\
\end{pmatrix};
\end{multline*}

\vspace{-0.5cm}

\begin{multline*}
R(R(e_{31})e_{31} + e_{31}R(e_{31})) \\
 = R(c_{13}e_{11} + c_{23}e_{21} + c_{11}e_{31}
 + c_{11}e_{31} + (c_{23}-a_{13})e_{32} + c_{13}e_{33} ) \\
 = 2c_{11}R(e_{31})
 + c_{13}\bigg(e_{12}+\frac{e_{23}}{2}\bigg)
 + c_{23}(-e_{11}+a_{13}e_{13})
 - 1/2(c_{23}-a_{13})(e_{11}+e_{22}).
\end{multline*}

For $e_{33}$, we have $c_{11}^2 = 2c_{11}^2$, i.\,e., $c_{11} = 0$.
For $e_{23}$, we have $-c_{13}/2 = c_{13}/2$ and so $c_{13} = 0$.
For $e_{11}$, we have $-(c_{23}-a_{13})/2 = -c_{23}-(c_{23}-a_{13})/2$,
it implies $c_{23} = 0$.

So, we have the RB-operator $R$ satisfying
\begin{gather*}
R(e_{13}) = R(e_{23}) = R(e_{33}) = 0,\quad
R(e_{12}) = \frac{1}{2}e_{13},\quad
R(e_{11}) = e_{12} + \frac{1}{2}e_{23},\quad
R(e_{22}) = \frac{1}{2}e_{23},\\
R(e_{21})
 = \begin{pmatrix}
-1 & 0 & a \\
0 & 0 & 0 \\
0 & 0 & 0 \\
\end{pmatrix},\quad
R(e_{32})
 = \begin{pmatrix}
-\frac{1}{2} & 0 & 0 \\
0 & -\frac{1}{2} & 0 \\
0 & 0 & 0 \\
\end{pmatrix},\quad
R(e_{31})
 = \begin{pmatrix}
0 & -a & 0 \\
-\frac{1}{2} & 0 & 0 \\
0 & 0 & 0 \\
\end{pmatrix},
\end{gather*}
where $a = a_{13}$.
Actually, \texttt{Singular} applied for~\eqref{R1=e12+e23-FirstForm} gives exactly these coefficients.

Conjugation with $\psi_{1,1/2}$~\eqref{psi} maps $R$ to the RB-operator $P$ on $M_3(F)$ defined as
\begin{gather*}
P(e_{13}) = P(e_{23}) = P(e_{33}) = 0,\quad
P(e_{12}) = e_{13},\quad
P(e_{11}) = e_{12} + e_{23},\quad
P(e_{22}) = e_{23},\\
P(e_{21}) = -e_{11} + qe_{13},\quad
P(e_{32}) = -(e_{11}+e_{22}),\quad
P(e_{31}) = -e_{21} - q e_{12},
\end{gather*}
where $q = 2a$.

\vfill \eject

Finally, conjugation with an automorphism $\psi$ of $M_3(F)$, where
\begin{gather*}
\psi(e_{12})=e_{12},\quad
\psi(e_{13})=e_{13},\quad
\psi(e_{23})=e_{23},\\
\psi(e_{11})=e_{11}-qe_{13},\quad
\psi(e_{22})=e_{22},\quad
\psi(e_{33})=e_{33}+qe_{13},\\
\psi(e_{21})=e_{21}-qe_{23},\quad
\psi(e_{31})=e_{31}+q(e_{11}-e_{33})-q^2e_{13},\quad
\psi(e_{32})= qe_{12} + e_{32},
\end{gather*}
maps $P$ to the RB-operator $O$ on $M_3(F)$ defined as follows,
\begin{gather*}
O(e_{13}) = O(e_{23}) = O(e_{33}) = 0,\quad
O(e_{12}) = e_{13},\quad
O(e_{11}) = e_{12} + e_{23},\quad
O(e_{22}) = e_{23},\\
O(e_{21}) = -e_{11},\quad
O(e_{32}) = -(e_{11}+e_{22}),\quad
O(e_{31}) = -e_{21}.
\end{gather*}
The RB-operator $O$ is defined by Example~1 for $n = 3$.
Here $\psi(X) = A^{-1}XA$ for $A = E - qe_{13}$.

{\bf Lemma 2}.
Let $R$ be an RB-operator of weight~0 on $M_3(\mathbb{C})$ satisfying $(R(1))^2\neq0$,
then up to conjugation with transpose and with an automorphism of $M_3(\mathbb{C})$ and up to scalar multiple, $R$~is defined by Example~1.

{\bf Remark 1}.
Note that we have described all RB-operators on $M_3(\mathbb{C})$ of weight zero
of the maximal rb-index (equal to~5).

\section{Case $(R(1))^2 = 0$ and $R(1)\neq0$}

Now we arrive at the case $(R(1))^2 = 0$ for nonzero $R(1)$.
So, we may assume that
\begin{equation}\label{UnitImage-2}
R(1) = e_{12},\quad R(e_{12}) = 0.
\end{equation}

An automorphism~$\varphi$ of $M_3(F)$ such that $\varphi(e_{12}) = e_{12}$ has the following form:
\begin{equation} \label{AutoGen}
\begin{gathered}
e_{11}	\to \begin{pmatrix}
1 & \alpha & \beta \\
0 & 0 & 0 \\
0 & 0 & 0 \\
\end{pmatrix}, \quad
e_{22} \to \begin{pmatrix}
0 & -\alpha-\beta\delta & 0 \\
0 & 1 & 0 \\
0 & \delta & 0 \\
\end{pmatrix}, \quad
e_{33} \to \begin{pmatrix}
0 & \beta\delta & -\beta \\
0 & 0 & 0 \\
0 & -\delta & 1 \\
\end{pmatrix}, \\
e_{12}\to e_{12}, \
e_{13}\to \lambda({-}\delta e_{12} {+} e_{13}), \
e_{32}\to \frac{1}{\lambda}( e_{32} {-} \beta e_{12} ), \
e_{31}\to
\frac{1}{\lambda}\!\!\begin{pmatrix}
{-}\beta & {-}\alpha\beta & {-}\beta^2 \\
0 & 0 & 0 \\
1 & \alpha & \beta \\
\end{pmatrix}\!\!,  \\
e_{21}\to \!\begin{pmatrix}
{-}\alpha{-}\beta\delta & {-}\alpha^2{-}\alpha\beta\delta & {-}\alpha\beta{-}\beta^2\delta \\
1 & \alpha & \beta \\
\delta & \alpha\delta & \beta\delta \\
\end{pmatrix}\!\!, \
e_{23}\to
\!\begin{pmatrix}
0 & \delta\lambda ( \alpha + \beta\delta ) & {-}\lambda ( \alpha + \beta\delta ) \\
0 & -\delta \lambda & \lambda \\
0 & -\delta^2 \lambda & \lambda\delta \\
\end{pmatrix}\!\!,
\end{gathered}
\end{equation}
where $\lambda\neq0$.
Below, we will refer on this automorphism~$\varphi$ as $\varphi(\alpha,\beta,\lambda,\delta)$
with all four parameters specified.
Note that $\varphi(X) = A^{-1}XA$ for 
$A = \begin{pmatrix}
1 & \alpha & \beta \\
0 & 1 & 0 \\
0 & -\delta\lambda & \lambda 
\end{pmatrix}$.

In what follows, we will also apply the following automorphism $\Theta_{12}$ of $M_3(F)$.
The $\Theta_{12}$ acts on matrix unities as $\Theta_{12}(e_{ij}) = e_{i'j'}$,
where $1' = 2$, $2' = 1$, and $3' = 3$.
Note that $\Theta_{12}\circ T$, where $T$ denotes the transpose of a matrix,
provides an antiautomorphism of $M_3(F)$ such that $(\Theta_{12}\circ T)(e_{12}) = e_{12}$.

Analogously to~\eqref{RxR1} and~\eqref{R1Rx}, we have
\begin{gather}
R(x)e_{12}
 = R(x)R(1) = R^2(x) + R(xe_{12}), \label{RxR1-2} \\
e_{12}R(x)
 = R(1)R(x) = R^2(x) + R(e_{12}x) \label{R1Rx-2}
\end{gather}
for any $x\in M_3(\mathbb{C})$.

For $x\in X = \{e_{11},e_{13},e_{22},e_{32},e_{33}\}$, by~\eqref{RxR1-2} and~\eqref{R1Rx-2} we get
$R(x)e_{12} = R^2(x) = e_{12}R(x)$, i.\,e. the following equality holds for $Y = (y_{ij}) = R(x)$:
\begin{equation}\label{R^2x-2}
R^2(x)
 = R(x)e_{12}
 = \begin{pmatrix}
0 & y_{11} & 0 \\
0 & y_{21} & 0 \\
0 & y_{31} & 0 \\
\end{pmatrix}
 = \begin{pmatrix}
y_{21} & y_{22} & y_{23} \\
0 & 0 & 0 \\
0 & 0 & 0 \\
\end{pmatrix}
 = e_{12}R(x).
\end{equation}
Thus, $R(x)  = \begin{pmatrix}
y_{11} & y_{12} & y_{13} \\
0 & y_{11} & 0 \\
0 & y_{32} & y_{33} \\
\end{pmatrix}$ for all $x\in X$.

Introduce matrices $A = R(e_{13})$, $B = R(e_{32})$, $C = R(e_{22}-e_{11})$ and $D = R(e_{11})$,
\begin{gather*}
R(e_{13})
 = \begin{pmatrix}
 a_{11} & a_{12} & a_{13} \\
 0 & a_{11} & 0 \\
 0 & a_{32} & a_{33} \\
 \end{pmatrix},\quad
R(e_{32})
 = \begin{pmatrix}
 b_{11} & b_{12} & b_{13} \\
 0 & b_{11} & 0 \\
 0 & b_{32} & b_{33} \\
 \end{pmatrix},\\
R(e_{22}-e_{11})
 = \begin{pmatrix}
 c_{11} & c_{12} & c_{13} \\
 0 & c_{11} & 0 \\
 0 & c_{32} & c_{33} \\
 \end{pmatrix},\quad
R(e_{11})
 = \begin{pmatrix}
 d_{11} & d_{12} & d_{13} \\
 0 & d_{11} & 0 \\
 0 & d_{32} & d_{33} \\
 \end{pmatrix}.
\end{gather*}

From~\eqref{RxR1-2} and~\eqref{R1Rx-2}, we have
$$
[R(e_{21}),e_{12}] = C, \quad
[R(e_{23}),e_{12}] = -A, \quad
[R(e_{31}),e_{12}] = B.
$$
For $T = R(e_{21})$, we have
$$
\begin{pmatrix}
0 & t_{11} & 0 \\
0 & t_{21} & 0 \\
0 & t_{31} & 0 \\
\end{pmatrix}
 - \begin{pmatrix}
 t_{21} & t_{22} & t_{23} \\
 0 & 0 & 0 \\
 0 & 0 & 0 \\
 \end{pmatrix}
 = \begin{pmatrix}
 c_{11} & c_{12} & c_{13} \\
 0 & c_{11} & 0 \\
 0 & c_{32} & c_{33} \\
 \end{pmatrix},
$$
so, $c_{11} = c_{33} = t_{21} = 0$ and
$$
R(e_{21})
 = \begin{pmatrix}
t_{22}+c_{12} & t_{12} & t_{13} \\
0 & t_{22} & -c_{13} \\
c_{32} & t_{32} & t_{33} \\
\end{pmatrix}
$$
Denote $F = R(e_{23})$ and $S = R(e_{31})$.
Analogously, we get $a_{11} = a_{33} = f_{21} = b_{11} = b_{33} = s_{21} = 0$ and so,
$$
R(e_{23}) = \begin{pmatrix}
f_{11} & f_{12} & f_{13} \\
0 & f_{11}+a_{12} & a_{13} \\
-a_{32} & f_{32} & f_{33} \\
\end{pmatrix}, \quad
R(e_{31}) = \begin{pmatrix}
s_{22}+b_{12} & s_{12} & s_{13} \\
0 & s_{22} & -b_{13} \\
b_{32} & s_{32} & s_{33} \\
\end{pmatrix}.
$$

From the equality
$$
R(e_{13})R(e_{32})
 = a_{13}b_{32}e_{12}
 = R(R(e_{13})e_{32} + e_{13}R(e_{32})) = 0
$$
and the symmetric equality for $R(e_{32})R(e_{13})$,
we get the relations
$$
a_{13}b_{32} = 0, \quad
a_{32}b_{13} = 0.
$$

If $a_{32}\neq0$, then $b_{13} = 0$. Since $R$ is nilpotent, then $b_{32} = 0$
and $R(e_{32}) = b_{12}e_{12}$. Otherwise, $R(e_{13}) = a_{12}e_{12}$.
Applying the conjugation with $\Theta_{12}\circ T$,
we may assume that $R(e_{13}) = ae_{12}$ (here $a = a_{12}$).
From the nilpotency of $R$, we deduce $b_{32} = 0$.

From
$$
a(s_{22}+b_{12})e_{12}
 = R(e_{31})R(e_{13})
 = R((s_{22}+b_{12})e_{13})
 + aR(e_{32}),
$$
we derive that $aR(e_{32}) =  0$.
Thus, either $R(e_{13}) = 0$ or $R(e_{32}) = 0$.
Again we may assume that $R(e_{13}) = 0$, i.\,e., $a = 0$.

Suppose that $c_{32}\neq0$. Since $\ker(R)$ is an $\Imm(R)$-module,
$e_{13}R(e_{21}),R(e_{21})e_{13}\in \ker(R)$. Thus,
$e_{11},e_{33}\in\ker(R)$. So, $R(e_{22}) = R(1) = e_{12}$,
and $c_{32} = 0$, a contradiction. We have proved that $c_{32} = 0$.

From the equality
\begin{multline*}
\begin{pmatrix}
d_{11}^2 & 2d_{11}d_{12}+d_{13}d_{32} & d_{13}(d_{11}+d_{33}) \\
0 & d_{11}^2 & 0 \\
0 & d_{32}(d_{11}+d_{33}) & d_{33}^2 \\
\end{pmatrix}
 = R(e_{11})R(e_{11}) \\
 = R(2d_{11}e_{11}+d_{12}e_{12}+d_{13}e_{13})
 = \begin{pmatrix}
 2d_{11}^2 & 2d_{11}d_{12}+d_{13} & 2d_{11}d_{13} \\
 0 & 2d_{11}^2 & 0 \\
 0 & 2d_{11}d_{32} & 2d_{11}d_{33} \\
 \end{pmatrix},
\end{multline*}
we conclude that $d_{11} = d_{33} = 0$.

From
\begin{multline*}
\begin{pmatrix}
0 & d_{12}f_{11}+d_{13}f_{32} & d_{13}f_{33} \\
0 & 0 & 0 \\
0 & d_{32}f_{11} & 0
\end{pmatrix}
 = \begin{pmatrix}
 0 & d_{12} & d_{13} \\
 0 & 0 & 0 \\
 0 & d_{32} & 0 \\
 \end{pmatrix}\begin{pmatrix}
 f_{11} & f_{12} & f_{13} \\
 0 & f_{11} & 0 \\
 0 & f_{32} & f_{33} \\
 \end{pmatrix}
 = R(e_{11})R(e_{23}) \\
 = R(d_{32}e_{33} + f_{11}e_{11})
 = f_{11}\begin{pmatrix}
 0 & d_{12} & d_{13} \\
 0 & 0 & 0 \\
 0 & d_{32} & 0 \\
 \end{pmatrix}
 + d_{32}\begin{pmatrix}
 0 & 1-2d_{12}-c_{12} & -2d_{13}-c_{13} \\
 0 & 0 & 0 \\
 0 & -2d_{32} & 0
\end{pmatrix},
\end{multline*}
we derive that $d_{32} = 0$.

Thus, we now have
\begin{equation}\label{21par-def}
\begin{gathered}
R(e_{12}) = R(e_{13}) = 0,\quad
R(e_{32}) = b_{12}e_{12} + b_{13}e_{13},\\
R(e_{22}-e_{11}) = c_{12}e_{12} + c_{13}e_{13}, \quad
R(e_{11}) = d_{12}e_{12} + d_{13}e_{13},\quad R(1) = e_{12}, \\
R(e_{23}) = \begin{pmatrix}
f_{11} & f_{12} & f_{13} \\
0 & f_{11} & 0 \\
0 & f_{32} & f_{33} \\
\end{pmatrix}, \\
R(e_{31}) = \begin{pmatrix}
s_{22}+b_{12} & s_{12} & s_{13} \\
0 & s_{22} & -b_{13} \\
0 & s_{32} & s_{33} \\
\end{pmatrix}, \quad
R(e_{21})
 = \begin{pmatrix}
t_{22}+c_{12} & t_{12} & t_{13} \\
0 & t_{22} & -c_{13} \\
0 & t_{32} & t_{33} \\
\end{pmatrix}.
\end{gathered}
\end{equation}

\texttt{Singular} allows to derive the following equations (it is not the full list of them):
\begin{gather}
\{b_{12},b_{13},c_{13},d_{13},s_{13},s_{22},s_{33},f_{13}+t_{22}-t_{33},2d_{12}+c_{12}-1\}\times\{f_{11},f_{33},f_{32}\} = 0, \\
\{b_{13},c_{13},d_{13},d_{12},t_{32},s_{22}-s_{33},t_{22}-t_{33},s_{32},f_{13}\}\times \{b_{12}-c_{13}\} = 0, \label{b12-c13} \\
\{b_{13},c_{13},d_{13}\}\times\{s_{12},s_{13},t_{12},t_{13},f_{12}+t_{32},3d_{12}+c_{12}-1+f_{13}\} = 0, \\
\{b_{13},c_{13},d_{13},b_{12}\}\times\{3d_{13}+2c_{13}+s_{22}-s_{33}\} = 0, \\
\{b_{13},c_{13},d_{13},s_{33}-s_{22}\}\times\{c_{12}+t_{22}-t_{33}+s_{32}\} = 0, \\
\{b_{13},c_{13},d_{13},2d_{12}+c_{12}-1,b_{12}\}\times\{f_{13}+s_{32}\} = 0, \\
\{f_{11}\}\times\{f_{33},t_{33}\} = 0.
\end{gather}

{\sc Case 1}:\,$b_{13}\neq0$.
Applying the conjugation with the\,automorphism\,$\varphi(0,0,1,-c_{13}/b_{13})$, we may assume that $c_{13} = 0$.
Hence, by~\eqref{b12-c13}, we get $b_{12} = 0$ too.
Moreover, $f_{11} = f_{33} = f_{32} = s_{12} = s_{13} = t_{12} = t_{13} = 0$. 
The computations with \texttt{Singular} also imply 
$$
f_{12}
 = s_{32} + c_{12}
 = t_{32}
 = t_{22} - t_{33}
 = 3d_{12}+2c_{12}-1
 = f_{13} - c_{12}
 = s_{22} - s_{33} + 3d_{13}
 = 0.
$$
Therefore, we have
\begin{gather*}
e_{12}, e_{13} \in \ker R,\quad
R(e_{32}) = b_{13}e_{13}, \quad
R(e_{22}-e_{11}) = c_{12}e_{12}, \\
R(e_{11}) = \frac{1-2c_{12}}{3}e_{12} + d_{13}e_{13},\quad R(1) = e_{12}, \quad
R(e_{23}) = c_{12}e_{13}, \\
R(e_{31}) = \begin{pmatrix}
s_{33}-3d_{13} & 0 & 0 \\
0 & s_{33}-3d_{13} & -b_{13} \\
0 & -c_{12} & s_{33} \\
\end{pmatrix}, \quad
R(e_{21})
 = \begin{pmatrix}
t_{33}+c_{12} & 0 & 0 \\
0 & t_{33} & 0 \\
0 & 0 & t_{33} \\
\end{pmatrix}.
\end{gather*}
Applying conjugation with~$\varphi(0,0,1/\sqrt{b_{13}},0)$, we may assume that $b_{13} = 1$.

Additionally, we have the six quadratic relations:
\begin{gather*}
(2c_{12}-1)(c_{12} + 1) = 0, \quad
d_{13}(c_{12}+1) = 0, \quad
6d_{13}^2+ 2c_{12}-1 = 0, \\
s_{33}(s_{33}-3d_{13}) = c_{12}+t_{33}, \quad
t_{33}(s_{33}-3d_{13}) = 0, \quad
t_{33}(t_{33}+c_{12}) = 0.
\end{gather*}

{\sc Case 1a}: $c_{12} = 1/2$.
Then $d_{13} = 0$. The rest system has the form
$$
s_{33}^2 = t_{33}+1/2, \quad
t_{33}s_{33} = 0, \quad
t_{33}(t_{33}+1/2) = 0.
$$

{\sc Case 1aa}: $t_{33} = -1/2$.
So, $s_{33} = 0$. After multiplication by~2 and application conjugation with~$\varphi(0,0,\sqrt{2},0)$ we have the RB-operator
\begin{equation}
\begin{gathered}
e_{11}, e_{12}, e_{13} \in \ker R, \quad
R(e_{32}) = e_{13}, \quad
R(e_{22}) = R(e_{33}) = e_{12}, \\
R(e_{23}) = e_{13}, \quad
R(e_{31}) = - e_{23} - e_{32}, \quad
R(e_{21}) = - e_{22} - e_{33}.
\end{gathered}
\end{equation}

{\sc Case 1ab}: $t_{33} = 0$.
So, $s_{33} = \pm 1/\sqrt{2}$. 
Applying conjugation with~$\varphi(0,0,1/s_{33},0)$, we get the operator
\begin{equation}
\begin{gathered}
e_{11}, e_{12}, e_{13} \in \ker R, \quad
R(e_{32}) = 2e_{13}, \quad
R(e_{2i}) = e_{1i}/2, \ i=1,2,3, \\
R(1) = e_{12}, \quad
R(e_{31}) = E - 2e_{23} - e_{32}/2.
\end{gathered}
\end{equation}

{\sc Case 1b}: $c_{12} = -1$. 
Hence, $d_{13} = \pm 1/\sqrt{2}$.

{\sc Case 1ba}: $t_{33} = 1$.
So, $s_{33} = 3d_{13}$. 
Applying conjugation with~$\varphi(0,0,d_{13},-d_{13})$, we get the RB-operator
\begin{equation} \label{eq4}
\begin{gathered}
e_{12}, e_{13}, e_{22}, e_{23} \in \ker R,\quad
R(e_{32}) = 2e_{13} - e_{12}, \quad
R(e_{11}) = e_{12}/2 + e_{13}, \\
R(e_{33}) = e_{12}/2 - e_{13}, \quad
R(e_{31}) = -2 e_{23} + e_{22} + 2e_{33}, \quad
R(e_{21}) = e_{22}/2 + e_{23}.
\end{gathered}
\end{equation}

{\sc Case 1bb}: $t_{33} = 0$. Hence, $s_{33} \in \{ d_{13}, 2d_{13} \}$.
Applying conjugation with an automorphism~$\varphi(0,0,d_{13},0)$,
we get the two RB-operators for $\kappa \in \{ 0, -1\}$:
\begin{equation}
\begin{gathered}
e_{12}, e_{13} \in \ker R,\quad
R(e_{32}) = 2 e_{13}, \quad
R(e_{22}-e_{11}) = -e_{12}, \quad
R(e_{11}) = e_{12} + e_{13}, \\
R(1) = e_{12}, \
R(e_{23}) = -e_{13}, \
R(e_{31}) = (\kappa{-}1)E - 2e_{23} + e _{32} + 3e_{33}, \
R(e_{21}) = -e_{11}.
\end{gathered}
\end{equation}

{\sc Case 2}: $b_{13} = 0$.
\texttt{Singular} implies that $c_{13} = d_{13} = 0$.
Moreover, we have the equations
$$
b_{12}\times \{ f_{11}, f_{32}, f_{33} \} = 0.
$$

Let us show that we may assume that
$f_{11} = f_{32} = f_{33} = 0$.
If $b_{12}\neq0$, then it follows immediately.
If $b_{12} = 0$, then one can find
among the equalities the following ones:
$$
\{ s_{22}, s_{33}, s_{13} \}
 \times \{ f_{11}, f_{32}, f_{33} \} = 0.
$$
Up to conjugation with $\Theta_{12}\circ T$,
we may assume that $f_{11} = f_{32} = f_{33} = 0$.
This application is correct, since here we have
$R(e_{13}) = R(e_{32}) = 0$.

Currently, we have
\begin{equation*}
\begin{gathered}
e_{12}, e_{13} \in \ker R,\
R(e_{32}) = b_{12}e_{12},\
R(e_{22}{-}e_{11}) = c_{12}e_{12}, \
R(e_{11}) {=} d_{12}e_{12},\
R(1) = e_{12}, \\
R(e_{23}) {=} f_{12}e_{12} + f_{13}e_{13}, \
R(e_{31}) {=} \begin{pmatrix}
s_{22}+b_{12} & s_{12} & s_{13} \\
0 & s_{22} & 0 \\
0 & s_{32} & s_{33} \\
\end{pmatrix}, \
R(e_{21}) {=} \begin{pmatrix}
t_{22}+c_{12} & t_{12} & t_{13} \\
0 & t_{22} & 0 \\
0 & t_{32} & t_{33} \\
\end{pmatrix}.
\end{gathered}
\end{equation*}

{\sc Case 2a}: $b_{12}\neq0$.
Then $d_{12} = f_{13} = s_{32} = t_{32} = 0$
and $s_{22} = s_{33}$, $t_{22} = t_{33}$.
Applying conjugation with~$\varphi(0,0,b_{12},-c_{12}/b_{12})$, 
we may assume that $b_{12} = 1$ and $c_{12} = 0$,
which implies $t_{13} = t_{22} = 0$ and $s_{12} = -s_{13}$.
Moreover, $t_{12} + f_{12}s_{13} = 0$ and $s_{22}\in \{ 0,-1\}$.
Conjugation with~$\varphi(0,s_{13},1,0)$ gives us the operators
\begin{equation}
\begin{gathered}
e_{11}, e_{12}, e_{13}, e_{21}, e_{22} \in \ker R,\quad
R(e_{32}) = e_{12},\quad
R(e_{33}) = e_{12}, \\
R(e_{23}) = f_{12}e_{12}, \quad
R(e_{31}) = \kappa E + e_{11}, \
\kappa \in \{ 0, -1 \}.
\end{gathered}
\end{equation}

{\sc Case 2b}: $b_{12} = 0$.
Among the relations, one has $d_{12}(d_{12}+c_{12}) = 0$.

{\sc Case 2ba}: $1 - 2d_{12} - c_{12} = 0$.
Then $(d_{12},c_{12}) = (0,1)$ or $(d_{12},c_{12}) = (1,-1)$.

{\sc Case 2baa}: $d_{12} = 0$ and $c_{12} = 1$. 
We have the equation~$f_{13}(f_{13}-1) = 0$.

{\sc Case 2baaa}: $f_{13} = 0$.
Applying conjugation with~$\varphi(0,0,1,f_{12})$,
we get that $f_{12} = 0$.
Then $s_{12} = s_{32} = t_{12} = t_{32} = 0$.

The rest equalities are
\begin{gather*}
s_{22} \times \{ s_{33}, t_{33} \} = 0, \quad
\{ s_{13}, s_{22}, s_{33}\} \times \{ t_{13}, t_{22}-t_{33}+1 \} = 0, \\
t_{13}(t_{22} - t_{33}) = 0, \quad
t_{33}( 2t_{22} - t_{33} + 1 ) = 0, \quad
t_{22}(t_{22}+1) = 0.
\end{gather*}
In particular, they imply that $t_{33}-t_{22} \in \{ 0, 1 \}$.

\vfill \eject

If $s_{22} \neq 0$, then we apply conjugation with~$\varphi(0,s_{13}/s_{22},s_{22},0)$
and obtain the operator
\begin{equation}
e_{11}, e_{12}, e_{13}, e_{23}, e_{32}, e_{33} \in \ker R, \
R(e_{22}) = e_{12}, \
R(e_{31}) = e_{11}+e_{22}, \
R(e_{21}) = - e_{22}.
\end{equation}

Let $s_{22} = 0$.
If $t_{22} = t_{33}$, we have the following operators after conjugation with an automorphism~$\varphi(0,t_{13},1,0)$:
\begin{equation} \label{eq8}
e_{11}, e_{12}, e_{13}, e_{23}, e_{31}, e_{32}, e_{33} {\in} \ker R, \
R(e_{22}) = e_{12}, \
R(e_{21}) = \kappa E + e_{11}, \ \kappa \in \{0,-1\}.
\end{equation}
Let $t_{33} = t_{22}+1$, hence, $t_{13} = 0$. 
If $s_{33} \neq 0$, we apply conjugation with~$\varphi(0,-s_{13}/s_{33},s_{33},0)$ and get the operators
\begin{equation}
\begin{gathered}
e_{11}, e_{12}, e_{13}, e_{23}, e_{32}, e_{33} \in \ker R, \quad
R(e_{22}) = e_{12}, \\
R(e_{31}) = e_{33}, \quad
R(e_{21}) = (\kappa+1)E - e_{22}, \ \kappa \in \{0,-1\}.
\end{gathered}
\end{equation}
Finally, let $s_{33} = 0$. Then we have the following operators (if $s_{13} \neq 0$, we apply conjugation with~$\varphi(0,0,\sqrt{s_{13}},0)$):
\begin{equation} \label{eq11}
\begin{gathered}
e_{11}, e_{12}, e_{13}, e_{23}, e_{32}, e_{33} \in \ker R, \quad
R(e_{22}) = e_{12}, \\
R(e_{31}) = - \mu e_{13}, \quad
R(e_{21}) = (\kappa+1)E - e_{22}, \ \kappa, \mu \in \{0,-1\}.
\end{gathered}
\end{equation}

{\sc Case 2baab}: $f_{13} \neq 0$, it means that $f_{13} = 1$. Then $f_{12} = 0$.

{\sc Case 2baaba}: $s_{13} = 0$.
Among the relations one has
$$
\{ s_{12}, t_{12}, t_{13} \} \times \{ s_{22}, s_{33} \} = 0.
$$

{\sc Case 2baabaa}: $(s_{12},t_{12},t_{13}) \neq (0,0,0)$.
Then $s_{22} = s_{33} = 0$. We have $s_{32}\in \{0, -1 \}$.

Suppose that $s_{32} = 0$. Then $s_{12} = 0$ and $t_{12} = t_{13}t_{32}$.
If $t_{13} = 0$, then $t_{12} = 0$, it is a contradiction to the conditions of Case 2baabaa.
Hence, $t_{13}\neq0$ and so $t_{33} = t_{22} + 1$.
Applying conjugation with~$\varphi(0,-t_{13},1,0)$,
we again come to the equality $t_{13} = 0$, which does not satisfy the conditions of the current subcase.

Now, consider the case $s_{32} = -1$.
If $s_{12} = 0$, then $t_{12} = t_{13} = 0$, again a contradiction. Thus, $s_{12}\neq0$.
We express $t_{32} = -\frac{t_{12}}{s_{12}}$, $t_{33} = t_{22} - \frac{t_{13}}{s_{12}}$.
The rest equations are
$$
t_{13} \times \{t_{12},t_{13}+s_{12}\} = 0, \quad
t_{22}(t_{22}+1) = 0.
$$
Applying conjugation with~$\varphi(0,0,1,-t_{12}/s_{12})$ and the equality $t_{12}t_{13} = 0$, we may assume that $t_{12} = 0$.
We get the following operators (depending on whether $t_{13} = 0$ or $t_{13} = -s_{12}$):
\begin{equation*}
\begin{gathered}
e_{11}, e_{12}, e_{13}, e_{32}, e_{33} \in \ker R,\quad
R(e_{22}) = e_{12}, \quad
R(e_{23}) = e_{13}, \\
R(e_{31}) = s_{12}e_{12} - e_{32}, \quad
R(e_{21}) = \kappa E + e_{11}, \ \kappa\in \{ 0,-1\};
\end{gathered}
\end{equation*}
\begin{equation*}
\begin{gathered}
e_{11}, e_{12}, e_{13}, e_{32}, e_{33} \in \ker R,\quad
R(e_{22}) = e_{12}, \quad
R(e_{23}) = e_{13}, \\
R(e_{31}) = s_{12}e_{12} - e_{32}, \quad
R(e_{21}) = (\kappa+1)E - e_{22} -s_{12}e_{13}, \ \kappa\in \{ 0,-1\}.
\end{gathered}
\end{equation*}
Applying conjugation with~$\varphi(0,s_{12},1,0)$, we reduce the case to the next one.

{\sc Case 2baabab}: $s_{12} = t_{12} = t_{13} = 0$.
Among the quadratic equalities, there are the following ones:
$(t_{22}-t_{33}+s_{32}+1)\times \{ s_{22},s_{33},t_{32} \} = 0$.

{\sc Case 2baababa}: $t_{33} = t_{22} + s_{32} + 1$.
The rest equations are
$$
s_{22} \times \{ s_{33}, t_{33} \} = 0, \quad
t_{22}^2 + t_{22} + t_{32}s_{22} = 0, \quad
t_{32}(s_{22}-s_{33})+s_{32}^2+s_{32} = 0.
$$

If $s_{22} = 0$ and $s_{33}\neq0$, then after conjugation
with~$\varphi(0,0,s_{33},-s_{32}/s_{33})$ we get the operators
\begin{equation}
\begin{gathered}
e_{11}, e_{12}, e_{13}, e_{32}, e_{33} \in \ker R,\quad
R(e_{22}) = e_{12}, \quad
R(e_{23}) = e_{13}, \\
R(e_{31}) = e_{33}, \quad
R(e_{21}) = (\kappa+1)E - e_{22}, \ \kappa \in \{0,-1\}.
\end{gathered}
\end{equation}

If $s_{22} = s_{33} = 0$, then we get the operators
(when $s_{32} = -1$, we apply conjugation with~$\varphi(0,0,1,t_{32})$;
when $s_{32} = 0$, we apply conjugation with~$\varphi(0,0,1,-t_{32})$)
\begin{equation} \label{eq13}
\begin{gathered}
e_{11}, e_{12}, e_{13}, e_{32}, e_{33} \in \ker R,\quad
R(e_{22}) = e_{12}, \quad
R(e_{23}) = e_{13}, \\
R(e_{31}) = -e_{32}, \quad
R(e_{21}) = \kappa E + e_{11}, \ \kappa \in \{ 0,-1\};
\end{gathered}
\end{equation}
\begin{equation} \label{eq14}
\begin{gathered}
e_{11}, e_{12}, e_{13}, e_{31}, e_{32}, e_{33} \in \ker R, \quad
R(e_{22}) = e_{12}, \quad
R(e_{23}) = e_{13}, \\
R(e_{21}) = (\kappa+1)E - e_{22}, \ \kappa \in \{ 0,-1\}.
\end{gathered}
\end{equation}

If $s_{22} \neq 0$, then $s_{33} = t_{33} = 0$ and we apply conjugation 
with~$\varphi(0,0,s_{22},(s_{32}+1)/s_{22})$ to write down the operator
\begin{equation}
\begin{gathered}
e_{11}, e_{12}, e_{13}, e_{32}, e_{33} \in \ker R,\
R(e_{2i}) = e_{1i}, \ i=1,2,3,\
R(e_{31}) = e_{11} + e_{22} - e_{32}.
\end{gathered}
\end{equation}

{\sc Case 2baababb}: $t_{33} \neq t_{22} + s_{32} + 1$.
Then $t_{32} = s_{22} = s_{33} = 0$, 
$t_{22}, s_{32} \in \{ 0, -1 \}$, and the remaining equality $t_{33}(2t_{22}-t_{33}+1) = 0$ implies that
$t_{33} = t_{22} - s_{32}$. Hence, we get the operator
\begin{equation} \label{eq16}
\begin{gathered}
e_{11}, e_{12}, e_{13}, e_{32}, e_{33} \in \ker R, \quad
R(e_{22}) = e_{12}, \quad
R(e_{23}) = e_{13}, \\
R(e_{31}) = \mu e_{32}, \quad
R(e_{21}) = \kappa E + e_{11} - \mu e_{33}, \ \kappa,\mu \in \{ 0, -1 \}.
\end{gathered}
\end{equation}

{\sc Case 2baabb}: $s_{13} \neq 0$.
Then $t_{33} = t_{22} + s_{32} + 1$.
Moreover,
$(t_{13}-s_{12})(s_{22}-s_{33}) = -s_{13}$.
Hence, $s_{22}-s_{33}\neq0$.
Conjugation with~$\varphi(0,s_{13}/(s_{22}-s_{33}),1,0)$
returns us to the case $s_{13} = 0$ (Case 2baaba).

{\sc Case 2bab}: $d_{12} = 1$ and $c_{12} = -1$.
Among the equations, we have
$f_{13}(f_{13}+1) = f_{12}(f_{13}+1) = 0$.

{\sc Case 2baba}: $f_{13} \neq 0$.
Then $f_{13} = -1$.
Applying~$\varphi(0,0,1,f_{12})$, we may assume that $f_{12} = 0$, which implies $s_{12} = t_{12} = t_{32} = 0$ and $s_{32} = 1$.

The rest equations are
\begin{gather*}
\{ s_{13}, s_{22}, s_{33} \} \times \{ t_{13}, t_{22} - t_{33} \} = 0, \quad
s_{22} \times \{ s_{33}, t_{33} \} = 0, \\
t_{33}(2t_{22} - t_{33} - 1) = 0, \quad
t_{13}(t_{22}-t_{33}-1) = 0, \quad
t_{22}(t_{22}-1) = 0.
\end{gather*}

{\sc Case 2babaa}: $t_{22} \neq t_{33}$.
Hence, $s_{13} = s_{22} = s_{33} = 0$ and $t_{33} = t_{22}-1$.
Thus, we get the operators
\begin{gather*}
e_{12}, e_{13}, e_{22}, e_{32}, e_{33} \in \ker R,\quad
R(e_{11}) = e_{12}, \quad
R(e_{23}) = -e_{13}, \\
R(e_{31}) = e_{32}, \quad
R(e_{21}) = \kappa E + e_{22} + t_{13}e_{13}, \ \kappa \in \{0,-1\}.
\end{gather*}
Applying the conjugation with~$\varphi(0,-t_{13},1,0)$, we may assume that $t_{13} = 0$.
The obtained operator is conjugate to~\eqref{eq13} with the help of~$\Theta_{12}$.

{\sc Case 2babab}: $t_{22} = t_{33}$.
Currently, we have
\begin{equation*}
\begin{gathered}
e_{12}, e_{13}, e_{22}, e_{32}, e_{33} \in \ker R,\quad
R(e_{11}) = e_{12}, \quad
R(e_{23}) = - e_{13}, \\
R(e_{31}) = \begin{pmatrix}
s_{22} & 0 & s_{13} \\
0 & s_{22} & 0 \\
0 & 1 & s_{33} \\
\end{pmatrix}, \quad
R(e_{21})
= \begin{pmatrix}
t_{22}-1 & 0 & 0 \\
0 & t_{22} & 0 \\
0 & 0 & t_{22} \\
\end{pmatrix}.
\end{gathered}
\end{equation*}

If $s_{22}\neq0$, then $s_{33} = t_{33} = 0$.
After conjugation with~$\varphi(0,s_{13}/s_{22},s_{22},0$) we get the operator
\begin{equation}
\begin{gathered}
e_{12}, e_{13}, e_{22}, e_{32}, e_{33} \in \ker R,\quad
R(e_{11}) = e_{12}, \quad
R(e_{23}) = -e_{13}, \\
R(e_{31}) = e_{11}+e_{22} + e_{32}, \quad
R(e_{21}) = -e_{11}.
\end{gathered}
\end{equation}

Let $s_{22} = 0$.
If $s_{33} \neq 0$, we apply conjugation with~$\varphi(0,-s_{13}/s_{33},s_{33},0)$ and get the operator
\begin{equation}
\begin{gathered}
e_{12}, e_{13}, e_{22}, e_{32}, e_{33} \in \ker R,\quad
R(e_{11}) = e_{12}, \quad
R(e_{23}) = -e_{13}, \\
R(e_{31}) = e_{32} + e_{33}, \quad
R(e_{21}) = (\kappa+1)E - e_{11}, \ \kappa \in \{ 0,-1\}.
\end{gathered}
\end{equation}

Let $s_{22} = s_{33} = 0$.
Depending on whether $s_{13}$ is zero, we write down the operators 
(if $s_{13}\neq0$, we conjugate with~$\varphi(0,0,\sqrt{s_{13}},0)$)
\begin{equation} \label{eq20}
\begin{gathered}
e_{12}, e_{13}, e_{22}, e_{32}, e_{33} \in \ker R,\quad
R(e_{11}) = e_{12}, \quad
R(e_{23}) = -e_{13}, \\
R(e_{31}) = e_{32} - \mu e_{13}, \quad
R(e_{21}) = (\kappa+1)E - e_{11}, \ \kappa,\mu \in \{ 0,-1\}.
\end{gathered}
\end{equation}
Note that the operator~\eqref{eq16} for $\mu = -1$
is conjugate by~$\Theta_{12}$ to~\eqref{eq20} with $\mu = 0$.

{\sc Case 2babb}: $f_{13} = 0$.
Then $f_{12} = 0$.
Let us show that we may assume that $s_{12} = 0$.
If $s_{32} = 1$, it follows by computations of~\texttt{Singular}. If $s_{32}\neq 1$, then we apply conjugation with~$\varphi(0,s_{12}/(1-s_{32}),1,0)$ and derive that $s_{12} = 0$.

Now, we want to show that we may consider only the case $t_{12} = 0$. Suppose that $t_{12}\neq0$.
Then $s_{13} = s_{22} = s_{33} = 0$, $t_{22} = t_{33}$, $s_{32} = 1$, $t_{12} = t_{13}t_{32}$.
Applying conjugation with~$\varphi(0,0,1,-t_{32})$,
we derive that $t_{12} = 0$.

{\sc Case 2babba}: $t_{13} = 0$.
Then $s_{13} = 0$.
If $s_{22} \neq0$, then $s_{33} = t_{33} = 0$, $t_{22}+s_{32}-1 = 0$,
$t_{32} = t_{22}(1-t_{22})/s_{22}$.
After conjugation with~$\varphi(0,0,s_{22},-t_{22}/s_{22}$) we get the operator
\begin{equation}
\begin{gathered}
e_{12}, e_{13}, e_{22}, e_{23}, e_{32}, e_{33} {\in} \ker R,\
R(e_{11}) {=} e_{12},\
R(e_{31}) = e_{11} {+} e_{22} {+} e_{32}, \
R(e_{21}) = -e_{11}.
\end{gathered}
\end{equation}

Let $s_{22} = 0$.
If $s_{33} = 0$, then conjugation with $\Theta_{1,2}$ converts this subcase to the one from Case 2baa.

If $s_{33}\neq0$, then $s_{32} = 1+t_{33}-t_{22}$.
The rest quadratic relations give
$t_{22} \in \{0,1\}$ and $t_{32} = s_{32}(s_{32}-1)/s_{33}$.
If $t_{22} = 0$, we apply conjugation with~$\varphi(0,0,s_{33},(1-s_{32})/s_{33})$ and get the operator
\begin{equation}
e_{12}, e_{13}, e_{22}, e_{23}, e_{32}, e_{33} \in \ker R,\
R(e_{11}) = e_{12}, \
R(e_{31}) = e_{32} + e_{33}, \
R(e_{21}) = -e_{11}.
\end{equation}
If $t_{22} = 1$, we apply conjugation with~$\varphi(0,0,s_{33},-s_{32}/s_{33})$ and get the operator
\begin{equation}
e_{12}, e_{13}, e_{22}, e_{23}, e_{32}, e_{33} \in \ker R,\quad
R(e_{11}) = e_{12}, \quad
R(e_{31}) = e_{33}, \quad
R(e_{21}) = e_{22}.
\end{equation}

{\sc Case 2babbb}: $t_{13} \neq 0$.
Then we derive
$$
s_{32} = 1, \quad t_{32} = 0, \quad t_{33} = t_{22}, \quad
s_{33} = s_{22} + \frac{s_{13}}{t_{13}}.
$$
Applying conjugation with~$\varphi(0,-t_{13},1,0)$,
we get $R(e_{21}) = t_{22}E - e_{11}$ and
$R(e_{31}) = s_{22}E + qe_{33} + e_{32}$.
Thus, we return to the case $t_{13} = 0$, which has been already studied in Case 2babba.

{\sc Case 2bb}: $1 - 2d_{12} - c_{12} \neq 0$.
Then $s_{13} = s_{22} = s_{33} = 0$ and $s_{32} = -f_{13}$.
Currently, we have
\begin{equation*}
\begin{gathered}
e_{12}, e_{13}, e_{32} \in \ker R,\quad
R(e_{22}-e_{11}) = c_{12}e_{12}, \quad
R(e_{11}) = d_{12}e_{12},\quad
R(1) = e_{12}, \\
R(e_{23}) = f_{12}e_{12} + f_{13}e_{13}, \quad
R(e_{31}) = \begin{pmatrix}
0 & s_{12} & 0 \\
0 & 0 & 0 \\
0 & -f_{13} & 0 \\
\end{pmatrix}, \quad
R(e_{21})
 = \begin{pmatrix}
t_{22}+c_{12} & t_{12} & t_{13} \\
0 & t_{22} & 0 \\
0 & t_{32} & t_{33} \\
\end{pmatrix}.
\end{gathered}
\end{equation*} 

We also have the equation $d_{12}(c_{12}+d_{12}) = 0$.
Up to conjugation with $\Theta_{12}\circ T$,
we may assume that $d_{12} = 0$.
By the conditions of the case, we have $c_{12} \neq 1$.
Hence, $t_{22} = t_{33}$. 
Moreover, $t_{13} = \dfrac{s_{12}(c_{12}-f_{13})}{c_{12}-1}$
and $t_{32} = \dfrac{f_{12}f_{13}}{c_{12}-1}$.
Now, we have
\begin{equation*}
\begin{gathered}
e_{11}, e_{12}, e_{13}, e_{32} \in \ker R, \quad
R(e_{22}) = c_{12}e_{12}, \quad
R(1) = e_{12}, \quad
R(e_{23}) = f_{12}e_{12} + f_{13}e_{13}, \\
R(e_{31}) = \begin{pmatrix}
0 & s_{12} & 0 \\
0 & 0 & 0 \\
0 & -f_{13} & 0 \\
\end{pmatrix}, \quad
R(e_{21}) = \begin{pmatrix}
t_{22}+c_{12} & t_{12} & \frac{s_{12}(c_{12}-f_{13})}{c_{12}-1} \\
0 & t_{22} & 0 \\
0 & \frac{f_{12}f_{13}}{c_{12}-1} & t_{22} \\
\end{pmatrix}.
\end{gathered}
\end{equation*} 
Applying conjugation with~$\varphi(0,0,1,f_{12}/(c_{12}-1))$, we get $t_{32} = 0$.

The rest equalities are
$$
f_{13} \times \{ f_{12}, t_{12}, f_{13}-c_{12} \} = 0, \quad
t_{22}(t_{22}+c_{12}) = 0, \quad
f_{12}t_{13} + c_{12}t_{12} = 0.
$$

{\sc Case 2bba}: $f_{13} \neq 0$.
Then $f_{12} = t_{12} = 0$ and $f_{13} = c_{12}\neq0$.

If $c_{12} \neq 1/2$, then we apply conjugation with~$\varphi(0,s_{12}/(2c_{12}-1),1,0)$
and obtain the operators
\begin{equation} \label{eq29}
\begin{gathered}
e_{11}, e_{12}, e_{13}, e_{32} \in \ker R, \quad
R(e_{22}) = c_{12}e_{12}, \quad
R(1) = e_{12}, \quad
R(e_{23}) = c_{12}e_{13}, \\
R(e_{31}) = -c_{12} e_{32}, \quad
R(e_{21}) = c_{12}(\kappa E + e_{11}), \ \kappa \in \{ 0, -1 \}.
\end{gathered}
\end{equation} 

Let $c_{12} = 1/2$.
Then we get after multiplication on~2 (if $s_{12} \neq 0$, we also apply conjugation with~$\varphi(0,0,s_{12},0)$) the operators
\begin{equation} \label{Q24}
\begin{gathered}
e_{11}, e_{12}, e_{13}, e_{32} \in \ker R, \quad
R(e_{22}) = R(e_{33}) = e_{12}, \quad
R(e_{23}) = e_{13}, \\
R(e_{31}) = - e_{32} - \mu e_{12}, \quad
R(e_{21}) = \kappa E + e_{11}, \ \kappa,\mu \in \{ 0, -1 \}.
\end{gathered}
\end{equation}

{\sc Case 2bbb}: $f_{13} = 0$.
Then $t_{12} = -\dfrac{s_{12}f_{12}}{c_{12}-1}$.
Currently, we have
\begin{equation*}
\begin{gathered}
e_{11}, e_{12}, e_{13}, e_{32} \in \ker R, \quad
R(e_{22}) = c_{12}e_{12}, \quad
R(1) = e_{12}, \quad
R(e_{23}) = f_{12}e_{12}, \\
R(e_{31}) = s_{12}e_{12}, \quad
R(e_{21}) = \begin{pmatrix}
t_{22}+c_{12} & -\frac{s_{12}f_{12}}{c_{12}-1} & \frac{s_{12}c_{12}}{c_{12}-1} \\
0 & t_{22} & 0 \\
0 & 0 & t_{22} \\
\end{pmatrix}.
\end{gathered}
\end{equation*}

If $c_{12}\neq 1/2$, we apply conjugation with~$\varphi(0,0,1,f_{12}/(2c_{12}-1))$ and get $f_{12} = 0$.
Further, we apply conjugation with~$\varphi(0,s_{12}/(c_{12}-1),1,0)$.
If $c_{12}\neq0$, then after dividing on~$c_{12}$ we arrive at
\begin{equation} \label{R31}
\begin{gathered}
e_{11}, e_{12}, e_{13}, e_{23}, e_{31}, e_{32} \in \ker R, \quad
R(e_{22}) = e_{12}, \\
R(e_{33}) = \alpha e_{12}, \quad
R(e_{21}) = \kappa E +e_{11}, \ \kappa \in \{ 0, -1\}, \ \alpha\neq0,\pm1.
\end{gathered}
\end{equation}
If $c_{12} = 0$, we write down the operator
\begin{equation} \label{dimImR=1}
e_{11}, e_{12}, e_{13}, e_{21}, e_{22}, e_{23}, e_{31}, e_{32} \in \ker R, \quad
R(e_{33}) = e_{12}.
\end{equation}

Let $c_{12} = 1/2$.
Conjugation with~$\varphi(0,-2s_{12},1,0)$ and multiplication on~2 gives us
\begin{equation*}
\begin{gathered}
e_{11}, e_{12}, e_{13}, e_{31}, e_{32} \in \ker R, \quad
R(e_{22}) = R(e_{33}) = e_{12}, \quad
R(e_{23}) = q e_{12}, \\
R(e_{21}) = \kappa E + e_{11}, \ \kappa \in \{ 0, -1 \}.
\end{gathered}
\end{equation*}
Hence, we get the following operators (if $q\neq 0$, we apply conjugation with~$\varphi(0,0,1/q,0)$):
\begin{equation}
\begin{gathered}
e_{11}, e_{12}, e_{13}, e_{23}, e_{31}, e_{32} \in \ker R, \\
R(e_{22}) = R(e_{33}) = e_{12}, \quad
R(e_{21}) = \kappa E + e_{11}, \ \kappa \in \{ 0, -1 \};
\end{gathered}
\end{equation}
\begin{equation}
\begin{gathered}
e_{11}, e_{12}, e_{13}, e_{31}, e_{32} \in \ker R, \quad
R(e_{22}) = R(e_{33}) = e_{12}, \quad
R(e_{23}) = e_{12}, \\
R(e_{21}) = \kappa E + e_{11}, \ \kappa \in \{ 0, -1 \}.
\end{gathered}
\end{equation}

We have considered all possible cases
and we are ready to formulate the main result of the work.
Below, we gather the operators~\eqref{eq13},~\eqref{eq29},~\eqref{Q24} with $\mu = 0$ as~(Q23) and~\eqref{eq8},~\eqref{R31},~\eqref{dimImR=1} as~(Q5). Moreover, we write only the nonzero action of an RB-operator on a matrix unity.

{\bf Theorem 3}.
Up to conjugation with automorphism or transpose and up to multiplication on a nonzero scalar,
one has the following nonzero Rota---Baxter operators of weight~0 on $M_3(F)$:
\begin{gather*}
(Q1)\
O(e_{12}) = e_{13},\quad
O(e_{11}) = e_{12} + e_{23},\quad
O(e_{22}) = e_{23}, \\
O(e_{21}) = -e_{11}, \quad
O(e_{32}) = -(e_{11}+e_{22}),\quad
O(e_{31}) = -e_{21}; \\
(Q2)\
R(e_{33}) = e_{12}; \\
(Q3)\
R(e_{32}) = R(e_{33}) = e_{12}, \quad
R(e_{23}) = \alpha e_{12}, \quad
R(e_{31}) = \kappa E + e_{11}; \\
(Q4)\
R(e_{22}) = R(e_{33}) = R(e_{23}) = e_{12}, \quad
R(e_{21}) = \kappa E + e_{11}; \\
(Q5)\
R(e_{22}) = e_{12}, \quad
R(e_{33}) = \beta e_{12}, \quad
R(e_{21}) = \kappa E +e_{11}; \\
(Q6)\
R(e_{22}) = e_{12}, \quad
R(e_{31}) = {-} \mu e_{13}, \quad
R(e_{21}) = (\kappa+1)E {-} e_{22}; \\
(Q7)\
R(e_{11}) = e_{12}, \quad
R(e_{31}) = e_{33}, \quad
R(e_{21}) = e_{22}; \\
(Q8)\
R(e_{22}) = e_{12}, \quad
R(e_{31}) = e_{33}, \quad
R(e_{21}) = (\kappa+1)E - e_{22}; \\
(Q9)\
R(e_{22}) = e_{12}, \quad
R(e_{23}) = e_{13}, \quad
R(e_{21}) = (\kappa+1)E - e_{22}; \\
(Q10)\
R(e_{22}) = e_{12}, \quad
R(e_{31}) = e_{11}+e_{22}, \quad
R(e_{21}) = - e_{22}; \\
(Q11)\
R(e_{11}) {=} e_{12}, \quad
R(e_{31}) = e_{11} {+} e_{22} {+} e_{32}, \quad
R(e_{21}) = -e_{11}; \\
(Q12)\
R(e_{11}) = e_{12}, \quad
R(e_{31}) = e_{32} + e_{33}, \quad
R(e_{21}) = -e_{11}; \\
(Q13)\
R(e_{22}) = e_{12}, \quad
R(e_{23}) = e_{13}, \quad
R(e_{21}) = \kappa E + e_{11}; \\ \allowdisplaybreaks
(Q14)\
R(e_{32}) {=} R(e_{23}) {=} e_{13}, \
R(e_{22}) {=} R(e_{33}) {=} e_{12}, \
R(e_{31}) {=} {-} e_{23} {-} e_{32}, \
R(e_{21}) {=} {-} e_{22} {-} e_{33}; \\
(Q15)\
R(e_{32}) = 2 e_{13}, \quad
R(e_{22}) = e_{13}, \quad
R(e_{11}) = e_{12} + e_{13}, \quad
R(1) = e_{12}, \\
R(e_{23}) = -e_{13}, \quad
R(e_{31}) = (\kappa{-}1)E - 2e_{23} + e _{32} + 3e_{33}, \quad
R(e_{21}) = -e_{11}; \\
\!(Q16)\
R(e_{2i}) = \frac{e_{1i}}2, \,i{=}1,2,3, \
R(e_{32}) = 2e_{13}, \
R(e_{33}) = \frac{e_{12}}2, \
R(e_{31}) = E {-} 2e_{23} {-} \frac{e_{32}}2; \\
(Q17)
R(e_{32}) {=} 2e_{13} {-} e_{12},\!
R(e_{i1}) {=} \frac{e_{i2}}2 {+} e_{i3}, i{=}1,\!2,
R(e_{33}) {=} \frac{e_{12}}2 {-} e_{13},\!
R(e_{31}) {=} {-}2 e_{23} {+} e_{22} {+} 2e_{33}; \\
(Q18)\
R(e_{2i}) = e_{1i}, \ i=1,2,3,\quad
R(e_{31}) = e_{11} + e_{22} - e_{32}; \\
(Q19)\
R(e_{2i}) {=} e_{1i}, \, i=2,3,\
R(e_{31}) {=} e_{33}, \quad
R(e_{21}) {=} (\kappa{+}1)E {-} e_{22}; \\
(Q20)\
R(e_{11}) = e_{12}, \quad
R(e_{2i}) = -e_{1i}, \ i=1,3,\quad
R(e_{31}) = e_{11} + e_{22} + e_{32}; \\
\!\!(Q21)\
R(e_{11}) = e_{12}, \quad
R(e_{23}) = -e_{13}, \quad
R(e_{31}) = e_{32} + e_{33}, \quad
R(e_{21}) = (\kappa+1)E - e_{11}; \\
(Q22)\
R(e_{11}) = e_{12}, \quad
R(e_{23}) = -e_{13}, \quad
R(e_{31}) = e_{32} - \mu e_{13}, \quad
R(e_{21}) = (\kappa+1)E - e_{11}; \\
\!\!(Q23)\
R(e_{2i}) = e_{1i}, \ i=2,3, \quad
R(e_{33}) = \gamma e_{12}, \quad
R(e_{31}) = - e_{32}, \quad
R(e_{21}) = \kappa E + e_{11}; \\
(Q24)\,
R(e_{22}) = R(e_{33}) = e_{12}, \quad
R(e_{23}) = e_{13}, \quad
R(e_{31}) = {-} e_{32} + e_{12}, \quad
R(e_{21}) = \kappa E {+} e_{11},
\end{gather*}
where $\kappa, \mu \in \{0,-1\}$, $\alpha,\beta,\gamma\in F$, $\beta\neq-1$, $\gamma\neq-1$.
Moreover, all listed above operators lie in different orbits of the set of RB-operators of weight~0 on $M_3(F)$ under
conjugation with $\Aut(M_3(F))$, transpose or multiplication on a nonzero scalar.

{\sc Proof}.
Let us denote the operator~$(Qi)$ with $\kappa = -1$ (where it is applicable) as~$(Qi)_{-1}$. The notation~$(Qi)_0$ will refer to the case $\kappa = 0$.

The case, when $R(1)^2\neq0$, was considered in the previous section; we get the only operator~(Q1).
Now, we study the case when $R(1)\neq0$ and $R(1)^2 = 0$.
Up to conjugation with the corresponding automorphism, we may assume that $R(1) = e_{12}$.

The operator~(Q2) is the only operator satisfying the property~\underline{$\dim(\Imm R) = 1$}.

The operators~(Q3)--(Q5) and (Q6) with $\mu = 0$ are exactly the ones with \underline{$\dim(\Imm R) = 2$}, let us call them group~I.
The operator~$(Q6)_{-1}$ with $\mu = 0$ is the only operator in the group~I such that $\Imm R$ contains an idempotent of rank~1 and $R^2 \neq 0$.
Further, the operator~$(Q6)_0$ with $\mu = 0$ is the only operator in the group~I such that $\Imm R$ contains an idempotent of rank~2 and $R^2 = 0$.
For the operator~(Q3), but not for (Q5) and (Q6), there exists a~preimage~$v$ of a one-sided unit of~$\Imm R$ such that $\Rad(\Imm(R))v = 0$ or $v\Rad(\Imm(R)) = 0$. Here $\Rad(A)$ denotes the Jacobson radical of an associative algebra. Since $\dim(\Imm R) < \infty$,
in all cases $\Rad(\Imm R)$ is the maximal nilpotent ideal of $\Imm R$. 
Also, a pair of operators from (Q3)--(Q5) with not the same $\kappa$ may not be conjugate, since the corresponding algebras $\Imm R$ contain idempotents of different rank.
Let us show that~$(Q3)$ with different values of~$\alpha$ are not pairwise conjugate. When $\alpha\neq0$, one may find $u,v,e$ such that $u,v$ are nilpotent matrices of rank~1, $e$ is a one-sided unit of $\Imm(R)$, and moreover, $(0\neq)R(v)\in \Span\{R(1)\}$, $R(u) = e$,
$$
\Imm(R)u = v\Imm(R) = ev = uv = 0, \quad
R(1)vu = e, \quad
R(1)v \neq 0.
$$
Then the coefficient of $R(v)$ at $R(1)$ is a constant, which is invariant under conjugation.
The operators~(Q4) and~(Q5) for different values of~$\beta$ are also not conjugate, since we may consider two orthogonal idempotents, preimages of a~basic element~$v$ of~$\rad(\Imm R)$, and find the ratio of their coefficients at~$v$. Up to conjugation, this ratio is fixed.
Finishing the case $\dim(\Imm R) = 2$,
let us separate~(Q4) and~(Q5) with $\beta = 1$.
For the first one, we may find a nilpotent matrix~$u$ of rank~1
and a one-sided unit~$e$ of~$\Imm R$
such that $\Rad(\Imm R)$ is linearly spanned by $R(u)$
and $eu = ue = 0$.

The operators (Q6) with $\mu = -1$ and (Q7)--(Q13) are exactly the operators with \underline{$\dim(\Imm R) = 3$}, let us call them as group~II.
The operators~(Q7) and~$(Q8)_0$ are the only ones from~II such that $R^2 = 0$ and $\dim(\Rad(\Imm R)) = 1$.
We distinguish them as follows: $\ker R$ contains a~non-degenerate matrix in the case of~$(Q8)_0$ and does not contain for~(Q7). The operator~$(Q13)_0$ is the only operator from~II satisfying the conditions $R^2 = 0$, $\dim(\Rad(\Imm R)) = 2$, and $\Imm(R)$ contains only idempotents of rank~1. The operators~$(Q6)_0$ with $\mu = -1$ and~$(Q9)_0$ are the only operators from~II such that $R^2 = 0$, $\dim(\Rad(\Imm R)) = 2$, and $\Imm(R)$ contains an idempotent of rank~2.
We separate this pair as follows: $\ker R$ contains a~non-degenerate matrix in the case of~$(Q6)_0$ with $\mu = -1$ and does not contain it for~$(Q9)_0$.

The operator~$(Q13)_{-1}$ is the only one such that $R^2 \neq 0$, $\dim(\Rad(\Imm R)) = 2$, and $\Imm R$ contains an idempotent of rank~2.
The operators~$(Q6)_{-1}$ with $\mu = -1$ and~$(Q9)_{-1}$ are the only operators from~II such that $R^2 \neq 0$, $\dim(\Rad(\Imm R)) = 2$, and $\Imm(R)$ contains an idempotent of rank~1.
Again, $\ker R$ contains a~non-degenerate matrix in the case of~$(Q6)_{-1}$ with $\mu = -1$ and does not contain for~$(Q9)_{-1}$.
The rest of the not-yet-considered operators from~II we split into two groups: $J_1 = \{ (Q10), (Q11) \}$ and $J_2 = \{ (Q12), (Q8)_{-1} \}$. Note that $\Imm R$ for operators from $J_2$ but not from $J_1$ contains an element~$v$ such that $v\cdot \Rad(\Imm R) = \Rad(\Imm R)\cdot v = 0$. Operators from both~$J_1$ and $J_2$ we distinguish by the property of whether $\ker R$ contains a~non-degenerate matrix or not.

Finally, the operators~(Q14)--(Q24) satisfy the condition~\underline{$\dim(\Imm R) = 4$}, let us collect these operators in the group~III.
The operator~$(Q22)_0$ is the unique operator from~III such that $R^2 = 0$, $\dim(\Rad(\Imm R)) = 3$, and $\Imm R$ contains an idempotent of rank~2.
The operator~$(Q22)_{-1}$ is the unique operator from~III such that $R^2 \neq 0$, $\dim(\Rad(\Imm R)) = 3$, and $\Imm R$ contains an idempotent of rank~1.
The operators~(Q14) and~$(Q15)_0$ are the only ones satisfying~$\dim(\Imm R\cap \ker R) = 2$.
In the case of~$(Q15)_0$ but not of~(Q14), there exists an idempotent~$e\in \Imm R$ of rank~1 such that $e$ acts on $\Rad(\Imm R)$ as a one-sided unit.

The operators~$(Q19)_0$ and~$(Q21)_0$ are the only ones from the rest of the operators satisfying the conditions $R^2 = 0$ and $\dim(\Rad(\Imm R)) = 2$.
These operators may not be conjugate, since their images are not (anti)isomorphic.

The operators~$(Q23)_0$ and $(Q24)_0$ are the only operators such that $R^2 = 0$, $\dim(\Rad(\Imm R)) = 3$,
and $\Imm R$ contains an idempotent of rank~1.
The operator~$(Q23)_0$ may be separated from $(Q24)_0$ as follows.
Let $v$ be a nonzero central element of $\Rad(\Imm R)$; we may find two orthogonal idempotents~$e,f$ of rank~1 such that
$R(e) = kv$ and $R(f) = lv$ for some $k,l\neq 0$. Then $k/l$ is a constant, which equals~$\gamma$ for $(Q23)_0$. By this reason, the operators~$(Q23)_0$ for different values of~$\gamma$ are not pairwise conjugate.
Let us separate~$(Q23)_0$ with~$\gamma = 1$ and~$(Q24)_0$.
For the operator~$(Q23)_0$ with~$\gamma = 1$ but not for $(Q24)_0$, we may find a nilpotent matrix~$v$ of rank~1 such that
$R(v) = u$ and $vu = uv = 0$ for a nonzero~$u\in \Rad(\Imm R)$ satisfying the equalities $eu = ue = 0$ for an idempotent~$e$ of~$\Imm R$. The operators~$(Q23)_{-1}$ and $(Q24)_{-1}$ are the only operators such that $R^2 \neq 0$, $\dim(\Rad(\Imm R)) = 3$,
and $\Imm R$ contains an idempotent of rank~2. Similarly to the previous arguments, we separate all these operators from each other.

The operators~$(Q15)_{-1}$,~(Q16),~(Q18),~(Q20),~$(Q21)_{-1}$ satisfy the conditions $R^2\neq 0$, $\dim(\Imm R\cap \ker R) = 3$,
$\dim(\Rad(\Imm R)) = 2$, and $\Imm R$ is isomorphic to $\Span\{ e_{11}, e_{12}$, $e_{13}, e_{22} \}$.
Further,~(Q16) and~(Q18), on the one hand, and~$(Q15)_{-1}$,~(Q20),~$(Q21)_{-1}$, on the other hand,
have non-isomorphic subalgebras~$\Imm R\cap \ker R$.
Now, for~(Q16) but not for~(Q18), there exist orthogonal idempotents~$e,f$ of rank~1 such that $R(e) = R(f) = R(1)/2$.
In the case of~(Q20), there exists $v$ such that $R(v)$ equals a one-sided unit of $\Imm R$ and $R(1)v = 0$.
For~$(Q21)_{-1}$ but not for~$(Q15)_{-1}$, one may find
$u,v$ such that $R(u),R(v)$ are orthogonal idempotents and $uR(1) = 0$. The same argument separates~(Q17) and~$(Q19)_{-1}$.
\hfill $\square$

{\bf Remark 2}.
Operators~(Q2), $(Q3)_0$--$(Q6)_0$, (Q7), $(Q8)_0$, $(Q9)_0$, and $(Q13)_0$ are exactly the nonzero ones satisfying $R^2 = 0$.
The operator~(Q1) is the only one such that $R^5 = 0$ and $R^4 \neq 0$. The rest operators satisfy~$R^3 = 0$ and $R^2 \neq 0$.
Note that there are no RB-operators~$R$ of weight~0 on~$M_3(F)$ such that $R^4 = 0 $ and $R^3\neq 0$.

{\bf Remark 3}.
Let us clarify how RB-operators of weight~0 on~$M_2(F)$ (described in Theorem~2) are incorporated into the list from~Theorem~3.
Denote
$$
A_0 = \Span\{ e_{11}, e_{12}, e_{21}, e_{22},e_{33} \}, \quad
A_1 = \Span\{ e_{13},e_{23},e_{31},e_{32}\}.
$$
Then $M_3(F) = A_0\dotplus A_1$ is a $\mathbb{Z}_2$-graded  algebra.
Hence, for any RB-operator~$R$ of weight~0 on~$M_3(F)$,
its projection~$P:= R|_{A_0}$ is an RB-operator of weight~0 on~$A_0$, 
provided that the projection of $\Imm R$ on~$A_1$ has a trivial product
and $R(R(x)y), R(yR(x)) \subset A_0$ for all $x,y \in A_0$.
Denoting $B = \Span \{ e_{11},e_{12},e_{21},e_{22}\} \cong M_2(F)$,
the operator~$P|_B$ as a~projection of an RB-operator on an ideal~$B$ of~$A_0$
is again an RB-operator.
One may check how this works for operators from Theorem~3.
In particular, the operators~$(Q5)_0$, $\beta = 0$,
and~$(Q6)_{-1}$, $\mu = 0$, are nothing more than Rota---Baxter operators~(L3) and~(L4)
from Theorem~2 extended on the whole~$M_3(F)$. 

{\bf Remark 4}.
The reason for the appearance of the two-valued coefficients~$\kappa$ and~$\mu$ in Theorem~3 is not clear.
Note that the operators~(L3) and~(L4) may be joined as follows:
$$
e_{11},e_{12}\in \ker R, \quad
R(e_{22}) = e_{12},\quad R(e_{21}) = \kappa E + e_{11}, \  \kappa \in \{ 0,-1\}.
$$

Finally, we consider the classification of V.V. Sokolov~\cite{Sokolov}
of all skew-symmetric Rota---Baxter operators of weight~0 on $M_3(\mathbb{C})$.

{\bf Theorem 4}~\cite{Sokolov}.
Up to conjugation, transpose and scalar multiple, all
nonzero skew-symmetric RB-operators on $M_3(\mathbb{C})$ are

(R1) $R(e_{31}) = e_{23}$, $R(e_{32}) = - e_{13}$;

(R2) $R(e_{11}) = - e_{21} - e_{32}$, $R(e_{12}) = e_{11} + e_{31}$,
$R(e_{13}) = e_{12} - e_{21}$, $R(e_{21}) = - e_{31}$, $R(e_{22}) = - e_{32}$,
$R(e_{23}) = e_{11} + e_{22}$;

(R3) $R(e_{23}) = e_{22}$, $R(e_{22}) = - e_{32}$;

(R4) $R(e_{13}) = e_{12} - e_{21}$, $R(e_{12}) = - R(e_{21}) = e_{31}$,
$R(e_{23}) = e_{11} + e_{22}$, $R(e_{11}) = R(e_{22}) = - e_{32}$;

(R5) $R(e_{13}) = e_{12}$, $R(e_{21}) = - e_{31}$, $R(e_{23}) = e_{11} + e_{22}$, $R(e_{11}) = R(e_{22}) = - e_{32}$;

(R6) $R(e_{33}) = e_{32}$, $R(e_{23}) = - e_{33}$, $R(e_{13}) = e_{11} + e_{12}$, $R(e_{11}) = R(e_{21}) = - e_{31}$;

(R7) $R(e_{23}) = -e_{11}-e_{33}$, $R(e_{11}) = R(e_{33}) = e_{32}$;

(R8) $R(e_{13}) = e_{12}$, $R(e_{21}) = -e_{31}$, $R(e_{33}) = e_{32}$, $R(e_{23}) = - e_{33}$.

Let us clarify to which operators from Theorem~4
all skew-symmetric RB-operators~$R$ of weight~0 on~$M_3(\mathbb{C})$
such that $R(1) \neq 0$ are conjugate.
Note that only (R1) does not satisfy the condition $R(1) \neq 0$.
Since the operator~(R2) satisfies~$R(1)^2 \neq 0$, it has to be conjugate to~(Q1).
The operator~(R3) up to sign and conjugation with~$\Theta_{1,3}$ coincides with~$(Q6)_{-1}$, $\mu = 0$.
The operator~(R4) up to sign and conjugation with~$\Theta_{1,3}$ and then with~$\varphi(0,0,i,0)$ coincides with~$(Q14)$.
The operator~(R5) up to sign and conjugation with~$\Theta_{1,3}$ coincides with~$(Q23)_{-1}$ with $\gamma = 1$.
Conjugation with the matrix~$E-e_{12}$ converts the operator~(R6) to the following one denoted by $P$:
\begin{gather*}
e_{12},e_{22},e_{31},e_{32}\in \ker P, \quad
e_{11}\to -e_{31}-e_{32}, \quad
e_{13}\to e_{11} + 2e_{12},\\
e_{23}\to e_{11} - e_{33} + 2e_{12},\quad
e_{21}\to -2e_{31}-2e_{32},\quad
e_{33}\to e_{32}.
\end{gather*}
Further, we apply conjugation with~$\Theta_{2,3}\circ T$,
then with~$\Theta_{1,2}\circ T$, and finally with $\varphi(0,0,-1,-1)$ to get up to sign~$(Q15)_0$.
The operator~(R7) up to conjugation with $\Theta_{1,2}\circ T\circ \Theta_{1,3}$ coincides with $(Q5)_{-1}$, $\beta = 1$.
The operator~(R8) up to conjugation with $\Theta_{1,3}$ coincides with~$(Q22)_{-1}$, $\mu = 0$.

\section*{Acknowledgements}

The research was carried out within the framework of 
the Sobolev Institute of Mathematics state contract (project FWNF-2022-0002).

\noindent Vsevolod Gubarev \\
Sobolev Institute of Mathematics \\
Acad. Koptyug ave. 4, 630090 Novosibirsk, Russia \\
Novosibirsk State University \\
Pirogova str. 1, 630090 Novosibirsk, Russia \\
e-mail: wsewolod89@gmail.com

\end{document}